\documentclass{llncs}

       \newcommand{\R}{\ensuremath{\sf R\hspace*{-1.0ex}\rule{0.15ex}%
       {1.5ex}\hspace*{0.9ex}}}

       \newcommand{\LC}{{\sf \hspace*{2ex}\rule{0.15ex}%
       {1.5ex}\hspace*{-1.5ex}\rule{1.5ex}{0.15ex}\hspace*{0.5ex}}}
       \newcommand{\RC}{{\sf \hspace*{0.7ex}\rule{0.15ex}%
       {1.5ex}\hspace*{-0.15ex}\rule{1.5ex}{0.15ex}\hspace*{0.5ex}}}
       \newcommand{\be}{\begin{equation}}
       \newcommand{\ee}{\end{equation}}
       \newcommand{\bv}[1]{\mathrm{\mathbf{#1}}}

\usepackage{graphicx}
\usepackage{pdfpages,epstopdf}
\DeclareGraphicsExtensions{.pdf}
\usepackage{epsfig}
\usepackage{verbatim}
\begin{document}
\title{Conic Sections and Meet Intersections in Geometric Algebra}
\author{Eckhard M.S. Hitzer}
\institute{Department of Physical Engineering, University of Fukui, Japan
           \email{hitzer@mech.fukui-u.ac.jp}}
\date{25 May. 2004}
\maketitle

\begin{abstract}
This paper first gives a brief overview over some interesting descriptions of conic sections, 
showing formulations in the three geometric algebras of Euclidean spaces, 
projective spaces, and the conformal model of Euclidean space. 
Second the conformal model descriptions of a subset of conic sections are listed
in parametrizations specific for the use in the main part of the paper. 
In the third main part the meets of lines and circles, and of spheres and planes are 
calculated for all cases of real and virtual intersections. 
In the discussion special attention is on the 
hyperbolic carriers of the virtual intersections.
\end{abstract}

\section{Introduction}

\subsection{Previous Work}

D. Hestenes used geometric algebra to give in his textbook 
New Foundations for Classical Mechanics~\cite{DH:NFCM} 
a range of descriptions of conic sections. The basic five ways
of construction there are:
\begin{itemize}
\item the semi-latus rectum formula
\item with polar angles (ellipse)
\item two coplanar circles (ellipse)
\item two non-coplanar circles (ellipse)
\item second order curves depending on three vectors
\end{itemize}
Animated and interactive online illustrations for all this can be found 
in~\cite{EH:ITMpres}.
Reference~\cite{EH:ITMpres} also treats plane conic sections defined via 
Pascal's Theorem by five general points in a plane in 
\begin{itemize}
\item the geometric algebra of the 2+1 dimensional projective plane 
\item the conformal geometric algebra of the 2+2 dimensional conformal model of the
Euclidean plane. 
\end{itemize}
This was inspired by Grassmann's treatment of plane conic sections in terms
of five general points in a plane~\cite{HG:A1A2e}. In both cases the meet operation is used in
an essential way. The resulting formulas are quadratic in each of the five
conformal points. 

By now it is also widely known that the conformal geometric algebra model of Euclidean 
space allows for direct linear product representations~\cite{DLL:CGeo,HLR:GS1} 
of the following subset of conics:
 Points,
 pairs of points,
 straight lines,
 circles,
 planes and
 spheres.
It is possible to find direct linear product representations with 5 constitutive 
points for general plane conics by introducing the geometric algebra of a six 
dimensional Euclidean vector space~\cite{CP:IntCon}.

Beyond this the meet~\cite{DH:CAGC} operation allows to e.g. 
generate a circle from the intersection of two spheres
or a sphere and a plane. The meet operation is well defined no matter whether two spheres
truly intersect each other (when the distance of the centers is less then the sum of the radii
but greater than their difference), 
but also when they don't (when the distance of the centers is greater than the sum of the radii
or less than their difference). 

The meet of two non-intersecting circles in a plane can be interpreted as a virtual point pair
with a distance that squares to a negative real number~\cite{LD:ITMga2003}. 
(If the circles intersect, the square is positive.)

This leads to the following set of questions:
\begin{itemize}
\item How does this virtual point pair depend on the locations of the centers? 
\item What virtual curve is generated if we continuously increase the center to center distances?
\item What is the dependence on the radii of the circles?
\item Does the meet of a straight line (a cirle with infinite radius) with a circle 
also lead to virtual point pairs and a virtual locus curve 
(depending on the distance of straight line and circle)?
\item How is the three dimensional situation of the meet of two spheres or a plane
and a sphere related to the two dimensional setting?
\end{itemize}
All these questions will be dealt with in this paper.

\section{Background}

\subsection{Clifford's Geometric Algebra}
\label{ssc:CGA}

Clifford's \textit{geometric algebra} $Cl(n-q,q) = \R{}_{\mathit{n-q,q}}$ of a real $n$-dimensional vector space 
$\R^{\mathit{n-q,q}}$
can be defined with four \textit{geometric product} axioms~\cite{GC:LAC} for a canonical vector basis, 
which satisfies
\begin{enumerate}
\item $ \bv{e}_k^2 = +1 \,\,(1 \leq k \leq n-q),\,\, \bv{e}_k^2 = -1 \,\,(n-q < k \leq n).$ 
\item The square of a vector $\bv{x}=x^k \bv{e}_k \,\, (1 \leq k \leq n)$ 
is given by the reduced quadratic form
$$
  \bv{x}^2 = (x^k \bv{e}_k)^2 = \sum_k (x^k)^2 e_k^2,
$$
which supposes 
$e_k e_l + e_l e_k = 0,\,\, k \neq l,\,\, 1 \leq k,l \leq n.$ 
\item Associativity: $(e_k e_l) e_m = e_k (e_l e_m), 1 \leq k,l,m \leq n.$
\item $\alpha e_k = e_k \alpha \mbox{ for all scalars } \alpha \in \R.$
\end{enumerate}
A geometric algebra is an example of a graded algebra with a basis of 
real scalars, vectors, bivectors, ... , $n$-vectors (pseudoscalars), i.e. the grades
range from $k=0$ to $k=n$. 
The grade $k$ elements form a 
$\left( \begin{array}{c} n \\ k \end{array}\right)$ 
dimensional $k$-vector space.
Each $k$-vector is in one-to-one correspondence with a $k$ dimensional subspace of 
$\R^{\mathit{n-q,q}}$. A general multivector $A$ of $\R{}_{\mathit{n-q,q}}$ is a sum of 
its grade $k$ parts
$$
  A = \sum_{k=0}^n \langle A \rangle_k. 
$$
The grade zero index is often dropped for brevity: $\langle A \rangle = \langle A \rangle_0$. 
Negative grade parts $k < 0$ or elements with grades $k>n$ do not exist, they are zero.
By way of grade selection a number of practically useful products of multivectors $A,B,C$ is derived from 
the geometric product:
\begin{enumerate}
\item The {scalar product}
\be
  A \ast B = \langle AB \rangle.
\ee
\item The {outer product} 
\be
  A \wedge B = \sum_{k,l=0}^n \langle \langle A \rangle_k  \langle B \rangle_l \rangle_{k+l}.
  \label{eq:op}
\ee
\item The left contraction
\be
  A \LC B = \sum_{k,l=0}^n \langle \langle A \rangle_k  \langle B \rangle_l \rangle_{l-k}
\ee
which can also be defined by 
\be
  (C \wedge A) \ast B = C \ast (A \LC B)
  \label{eq:lcos}
\ee 
for all  $C \in \R{}_{\mathit{n-q,q}}.$
\item The right contraction
\be
  A \RC B = \sum_{k,l=0}^n \langle \langle A \rangle_k  \langle B \rangle_l \rangle_{k-l}
  = (\tilde{B}\LC \tilde{A})\tilde{\,}
  \label{eq:rcos}
\ee
or defined by 
\be
  A \ast (B \wedge C) = (A \RC B) \ast C
\ee  
for all $C \in \R{}_{\mathit{n-q,q}}.$ The tilde sign $\tilde{A}$ indicates the
 reverse order of all elementary vector products in every grade component $\langle A \rangle_k$.
\item Hestenes and Sobczyk's~\cite{DH:CAGC} inner product generalization
\begin{eqnarray}
  \langle A \rangle_k  \cdot \langle B \rangle_l 
  & = & \langle \langle A \rangle_k \langle B \rangle_l \rangle_{|k-l|} \,\, (k \neq 0, l \neq 0),
  \label{eq:HSip}
  \\
  \langle A \rangle_k  \cdot \langle B \rangle_l 
  & = & 0 \,\, (k = 0 \mbox{ or } l = 0).
  \label{eq:HSip0}
\end{eqnarray}
\end{enumerate}
The scalar and the outer product (already introduced by H. Grassmann) are well accepted. 
There is some debate about the use of the left and right contractions on one hand or
Hestenes and Sobczyk's "minimal" definition on the other hand as the preferred generalizations of the 
inner product of vectors~\cite{LD:IPofGA}. For many practical purposes (\ref{eq:HSip}) 
and (\ref{eq:HSip0}) are completely sufficient. But the exception for grade zero factors
(\ref{eq:HSip0}) needs always to be taken into consideration when deriving formulas involving
the inner product. Hestenes and Sobczyk's book~\cite{DH:CAGC} shows this in a number of places.
The special consideration for grade zero factors (\ref{eq:HSip0}) also becomes necessary in
software implementations. Beyond this, (\ref{eq:lcos}) and (\ref{eq:rcos}) show how the salar and 
the outer product already fully imply left and right contractions. 
It is therefore infact possible to
begin with a Grassmann algebra, introduce a scalar product for vectors, induce the (left or right) 
contraction and thereby define the geometric product, 
which generates the Clifford geometric algebra. 
It is also possible to give direct definitions of the (left or right) 
contraction~\cite{PL:CAS}.

\subsection{Conformal Model of Euclidean Space}

Euclidean vectors are given in an orthonormal basis $\{\bv{e}_1,\bv{e}_2,\bv{e}_3 \}$ of
$\R^{3,0}=\R^3$ as
\be
 \bv{p}=p_1\bv{e}_1+p_2\bv{e}_2+p_3\bv{e}_3, \,\,\, \bv{p}^2=p^2  \enspace.
\ee
One-to-one corresponding conformal points in the 3+2 dimensions of $\R^{4,1}$ are given as
\be
  P=\bv{p}+\frac{1}{2}p^2\bv{n}+\bv{\bar{n}},  \,\,\, 
  P^2=\bv{n}^2=\bv{\bar{n}}^2=0,  \,\,\,
  P \ast \bv{\bar{n}} =
  \bv{n} \ast \bv{\bar{n}} = -1
\ee
with the special conformal points of
$\bv{n}$ infinity,
$\bv{\bar{n}}$ origin.
This is an extension of the Euclidean space similar to the projective model
of Euclidean space. But in the conformal model extra dimensions are introduced 
both for origin and infinity. $P^2=0$ shows that the conformal model first restricts
$\R^{\textstyle 4,1}$ 
to a four dimensional null cone (similar to a light cone
in special relativity) and second the normalization condition 
$ P \ast \bv{\bar{n}} =-1$ 
further intersects this cone with a hyperplane.

We define the Minkowski plane pseudoscalar (bivector) as
\be
 N=\bv{n} \wedge \bv{\bar{n}}, \,\,\, N^2=1 \enspace.
\ee

By joining conformal points with the outer product (\ref{eq:op}) we can generate the subset of conics mentioned above: 
pairs of points, straight lines, circles, planes and 
spheres~\cite{DLL:CGeo,GS:Wachter,HLR:GS1,EH:KWA,HU:JGAP}.
Detailed formulas to be used in the rest of the paper are given in the 
following subsections~\cite{EH:ego}.

\begin{figure}
\begin{center}
  \includegraphics[width=0.4\textwidth]{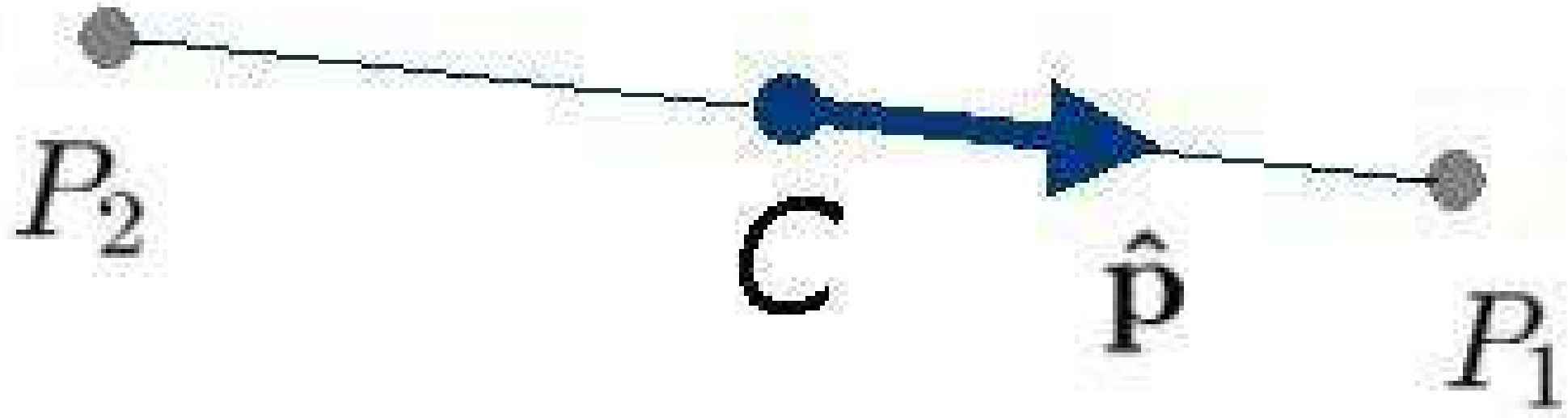}
  \caption[Conformal Pointpair]{
     Pair of intersection points $P_1, P_2$ with 
     distance $2r$, midpoint $C$ 
     and unit direction vector $\bv{\hat{p}}$ of the connecting line segment.
   }
  \label{fig:cpp}
\end{center}
\end{figure}

\subsection{Point Pairs} 

\be
   P_1\wedge P_2  
      =  \mathbf{p}_1\wedge\mathbf{p}_2
                     +\frac{1}{2}(p_2^2\mathbf{p}_1 - p_1^2\mathbf{p}_2)\mathbf{n}
                     -(\mathbf{p}_2-\mathbf{p}_1)\bar{\mathbf{n}}
                     +\frac{1}{2}(p_1^2- p_2^2)N
\ee
\be
   = \ldots = 2r \{\bv{\hat{p}}\wedge \bv{c} 
       + \frac{1}{2}[(c^2 + r^2)\bv{\hat{p}} -2\bv{c}\ast\bv{\hat{p}}\,\,\bv{c} ]\mathbf{n}
       + \bv{\hat{p}}\bv{\bar{n}}+\bv{c}\ast\bv{\hat{p}}N \} \enspace,
   \label{eq:PP}
\ee
with distance $2r$, unit direction of the line segment $\bv{\hat{p}}$, 
and midpoint $\bv{c}$ (comp. Fig.~\ref{fig:cpp}):
\be 
     2r =\, \mid\mathbf{p}_1-\mathbf{p}_2 \mid,     \,\,\, 
     \bv{\hat{p}}=\frac{\mathbf{p}_1-\mathbf{p}_2}{2r},  \,\,\,
     \bv{c}=\frac{\mathbf{p}_1+\mathbf{p}_2}{2} \enspace.
  \label{eq:PPrpc}
\ee
For $\bv{\hat{p}}\wedge \bv{c} =0 \,\,(\bv{\hat{p}} \,\| \,\bv{c})$ we get
\be
 P_1\wedge P_2=2r \{C-\frac{1}{2}\,r^2 \bv{n}\} \bv{\hat{p}} N
\ee
for $\bv{\hat{p}}\ast \bv{c}=0 \,\, (\bv{\hat{p}}\perp \bv{c})$ we get
\be
  P_1\wedge P_2=-2r \{C+\frac{1}{2}\,r^2 \bv{n}\}\bv{\hat{p}}
\ee
with conformal midpoint
\be
 C=\bv{c}+\frac{1}{2}c^2\bv{n}+\bv{\bar{n}} \enspace.
 \label{eq:cfcent}
\ee

\subsection{Straight Lines}

Using the same definitions the straight line through $P_1$ and $P_2$ is given by
\be
 P_1 \wedge P_2 \wedge \bv{n} 
   = 2r \bv{\hat{p}}\wedge C \wedge \bv{n}
   = 2r \{ \bv{\hat{p}}\wedge \bv{c}\,\,\bv{n} - \bv{\hat{p}}N \}
\ee

\subsection{Circles} 

$$
 P_1 \wedge P_2 \wedge P_3
   = \alpha \{\bv{c}\wedge I_c 
              +\frac{1}{2}[(c^2+r^2)I_c -2\bv{c} (\bv{c} \LC I_c)]\bv{n}  
              + I_c \bv{\bar{n}}   
              - (\bv{c} \LC I_c) N \} 
$$
\be
  = \alpha (C+\frac{1}{2}r^2\bv{n})\wedge \{I_c+\bv{n}(C\LC I_c)  \}
  \label{eq:cfcirc}
\ee
describes a circle through the three points $P_1$, $P_2$ and $P_3$ 
with center $\bv{c}$, radius $r$, and circle plane bivector
\be
  I_c = \frac{(\bv{p}_1-\bv{p}_2)\wedge(\bv{p}_2-\bv{p}_3)}{\alpha} \enspace,
\ee
where the scalar  $\alpha>0$ is chosen such that
\be
   I_c^2 =-1 \enspace.
\ee
For the inner product we use the left contraction $\LC$ 
in (\ref{eq:cfcirc}) as discussed in section \ref{ssc:CGA}.
For $\bv{c}\wedge I_c=0$ (origin $\bv{\bar{n}}$ in circle plane) and the
conformal center (\ref{eq:cfcent}) we get
\be
  P_1 \wedge P_2 \wedge P_3 = -\alpha\{C-\frac{1}{2}r^2\bv{n} \}I_c N.
\ee

\subsection{Planes} 

Using the same $\alpha$, $C$ and $I_c$ as for the circle
\be
  P_1 \wedge P_2 \wedge P_3 \wedge \bv{n}
    = \alpha \,\,C\wedge I_c \wedge \bv{n}
    = \alpha \{ \bv{c}\wedge I_c\bv{n} - I_c N  \}
  \label{eq:cplane}  
\ee
defines a plane through $P_1$, $P_2$, $P_3$ and infinity.
For $\bv{c}\wedge I_c=0$ (origin $\bv{\bar{n}}$ in plane) we get
\be
  P_1 \wedge P_2 \wedge P_3 \wedge \bv{n}
    = -\alpha  I_c N.
\ee

\subsection{Spheres} 

\be 
  P_1 \wedge P_2 \wedge P_3 \wedge P_4
    = \beta(C-\frac{1}{2}r^2\bv{n})I N
\ee
defines a sphere through $P_1$, $P_2$, $P_3$ and $P_4$ with radius $r$, 
conformal center $C$, unit volume trivector $I=\bv{e}_1\bv{e}_2\bv{e}_3$, 
and scalar 
\be
  \beta= (\bv{p}_1-\bv{p}_2)\wedge(\bv{p}_2-\bv{p}_3)\wedge(\bv{p}_3-\bv{p}_4)\,\,I^{-1}
  \enspace.
\ee

\section{Full Meet of Two Circles in One Plane}
\label{sc:2C}

\begin{figure}
\begin{center}
\includegraphics[width=1.0\textwidth]{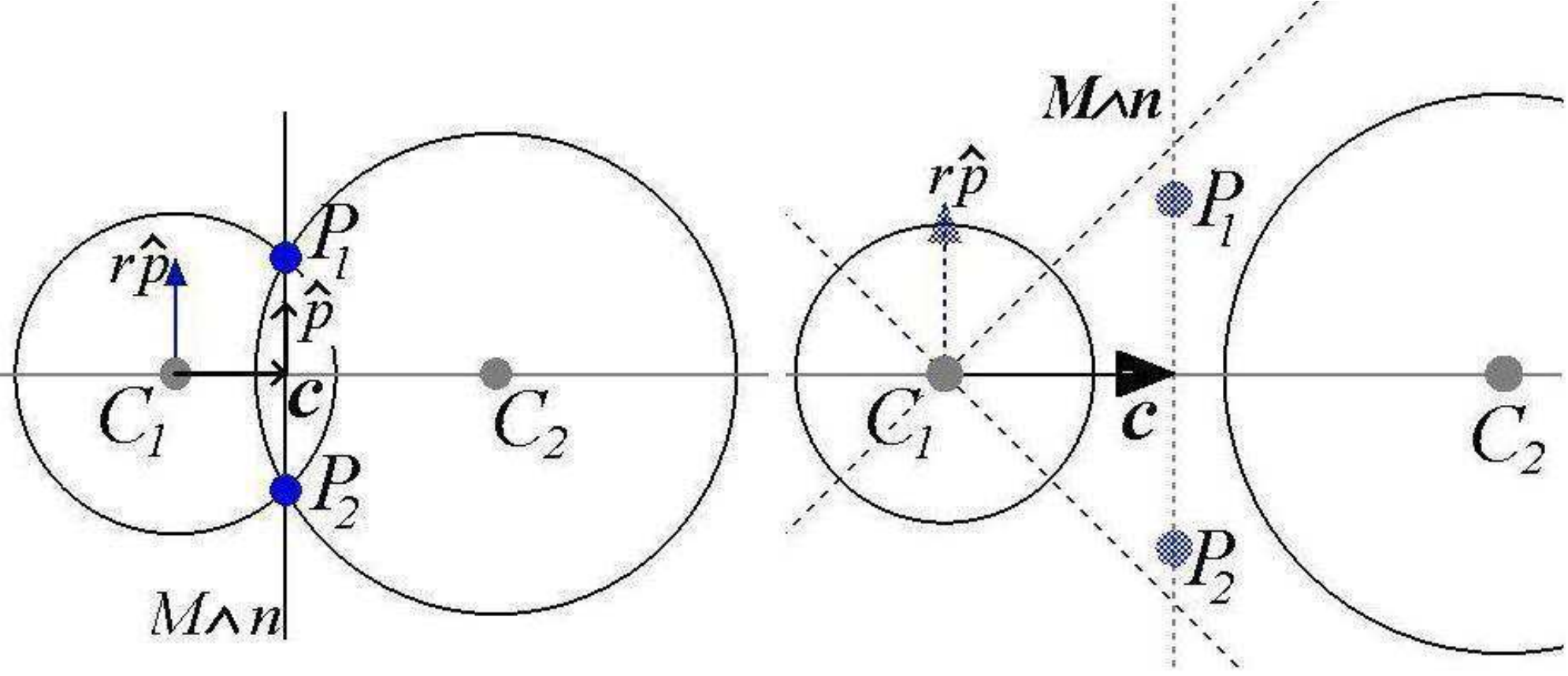}
  \caption[Circle circle intersection]{
  Two intersecting circles with centers $C_1,C_2$, intersecting in
  points $P_1,P_2$ at distance $2r$, with midpoint $\bv{c}$ 
  and unit direction vector $\bv{\hat{p}}$ of the connecting line segment.
  Left side: Real intersection ($r_2<d<r_1+r_2$), right side: Virtual intersection ($d>r_1+r_2$).}
  \label{fig:2c}
\end{center}
\end{figure}

\begin{figure}
\begin{center}
\resizebox{\textwidth}{!}{\includegraphics{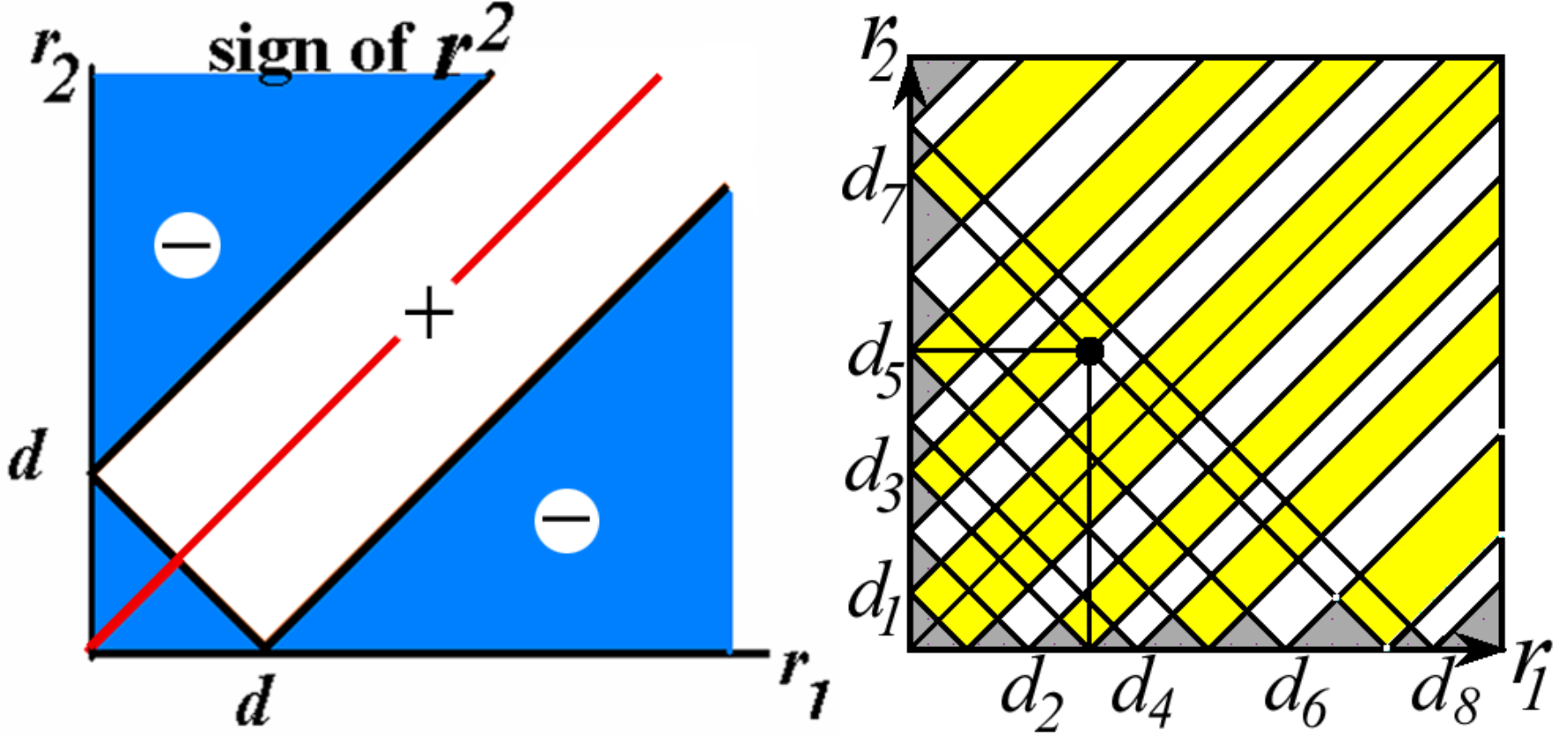}}
  \caption[r square sign]{
  Left side: Positive (+) and negative ($-$) signs of $r^2$ 
  depending on the radii of the two circles and the circle center distance $d$.
  $r=0$ (0) on the border lines (shape of an open U tilted in the $r_1=r_2$ direction), 
  which separate ($-$) and (+) regions.
  Right side: Black dot for case $r_1<r_2$. 
  Eight different values of $d$ are indicated, showing the tilted U-shaped 
  (0) border lines between the ($-$) and (+) regions. 
  $d_1<r_2-r_1$ ($-$),
  $d_2=r_2-r_1$ (0),
  $d_3=r_1$ (+),
  $r_2-r_1<d_4<r_2$ (+),
  $d_5=r_2$ (+),
  $r_2<d_6<r_1+r_2$ (+),
  $d_7=r_1+r_2$ (0) and
  $d_8>r_1+r_2$ ($-$).}
  \label{fig:r2s}
\end{center}
\end{figure}

The meet of two circles (comp. Fig.~\ref{fig:2c})
\be
  V_1 = (C_1-\frac{1}{2}r_1^2\bv{n})I_c N, \,\,\,\,\,\,
  V_2 = (C_2-\frac{1}{2}r_2^2\bv{n})I_c N \enspace,
\ee
with conformal centers $C_1, C_2$, radii $r_1, r_2$, in one plane $I_c$ 
(containing the origin $\bv{\bar{n}}$), and join four-vector $J = I_c N$
is\footnote{The dots ($\ldots$) in 
(\ref{eq:PP}), (\ref{eq:2cmeet}), 
(\ref{eq:clmeet}) 
and (\ref{eq:2smeet})
indicate nontrivial intermediate algebraic calculations whose details are
omitted here because of lack of space.} 
\be
  M = (V_1\LC J^{-1})\LC V_2
\ee
\be
  =\frac{1}{2}[(\bv{c}_1^2-r_1^2)-(\bv{c}_2^2-r_2^2)]I_c
   +\frac{1}{2}[(\bv{c}_2^2-r_2^2)\bv{c}_1 - (\bv{c}_1^2-r_1^2)\bv{c}_2 ]I_c \bv{n}
   +(\bv{c}_2-\bv{c}_1)I_c\bv{\bar{n}}
   +\bv{c}_1\wedge \bv{c}_2 I_c N
\ee
\be
  = \ldots = \frac{d}{2r} P_1\wedge P_2
  \label{eq:2cmeet}
\ee
with 
\be
  d=\mid \bv{c}_2-\bv{c}_1\mid \enspace,
  \label{eq:2cd}
\ee
\be
  r^2 = \frac{M^2}{(M\wedge \bv{n})^2}
      = d^2 \{\frac{r_1^2r_2^2}{d^4} - \frac{1}{4}(1-\frac{r_1^2}{d^2}-\frac{r_2^2}{d^2})^2 \}
\ee
and [like in (\ref{eq:PPrpc})]
\be
  \bv{p}_1=\bv{c}+r\bv{\hat{p}}, \,\,\,
  \bv{p}_2=\bv{c}-r\bv{\hat{p}} \enspace,
  \label{eq:2cp12}
\ee
\be
  \bv{c}=\bv{c}_1
         +\frac{1}{2}(1+\frac{r_1^2-r_2^2}{d^2})(\bv{c}_2-\bv{c}_1), \,\,\,\,\,\,
  \bv{\hat{p}}=\frac{\bv{c}_2-\bv{c}_1}{d}I_c \enspace.
\ee
We further get {\em independent} of $r^2$
\be
  M\wedge \bv{n} = d \bv{\hat{p}}\wedge C \wedge \bv{n} \enspace,
  \label{eq:2cMwn}
\ee
which is in general a
straight line through $P_1, P_2$, 
and in particular for $r^2=0$ $(M^2=0, \,\, \bv{p}_1=\bv{p}_2=\bv{c})$ 
the tangent line at the intersection point.\footnote{
We actually have 
$\lim_{r \rightarrow 0}M=\,\,d\{\bv{\hat{p}}+C\ast \bv{\hat{p}}\,\,\bv{n} \}\wedge C $,
which can be interpreted~\cite{LD:ITMga2003} as the tangent vector of the two
tangent circles, located at the point $C$ of tangency.
}
Note that $r^2$ may become {\em negative} (depending on $r_1$ and $r_2$, 
for details compare Fig.~\ref{fig:r2s}).
%
%
%
%
The vector from the first circle center $\bv{c}_1$ to the middle $\bv{c}$ of the point pair is
\be
  \bv{c}-\bv{c}_1 = \frac{1}{2}(1+\frac{r_1^2-r_2^2}{d^2})(\bv{c}_2-\bv{c}_1) \enspace,
\ee
with (oriented) length
\be
  d_1=\frac{1}{2}(d+\frac{r_1^2-r_2^2}{d})
\ee
This length $d_1$, half the intersection point pair distance $r$ 
and the circle radius $r_1$ are related by
\be
  r^2 + d_1^2 = r_1^2 \enspace.
  \label{eq:cchyp}
\ee

\begin{figure}
\begin{center}
\resizebox{\textwidth}{!}{\includegraphics{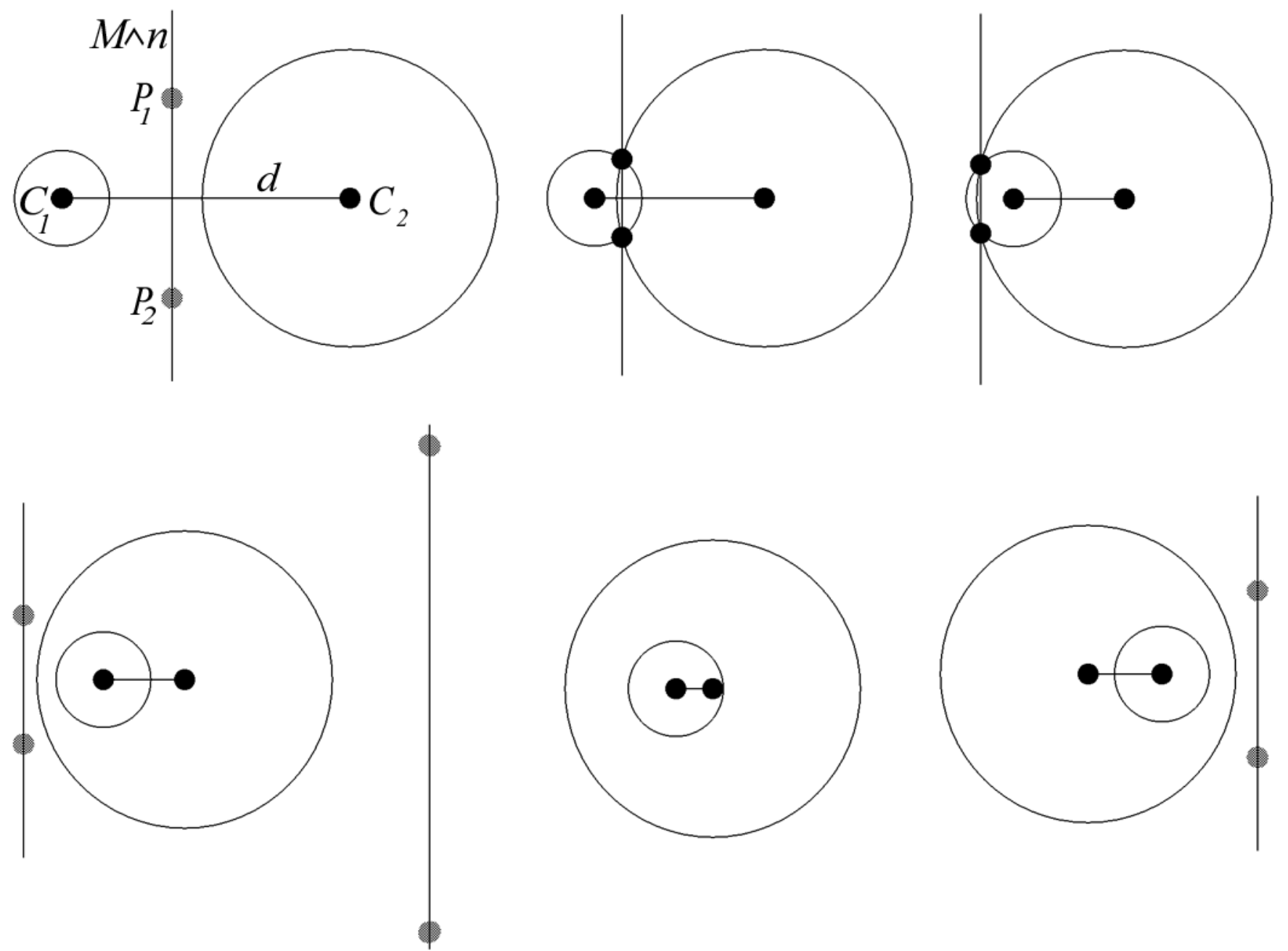}}
  \caption[Circle circle intersection series]{
  Real and virtual intersection points (\ref{eq:2cp12}) and
  vertical carrier line $M \wedge \bv{n}$ of (\ref{eq:2cMwn})
  for two circles with radii $r_1<r_2$ and centers $C_1,C_2$ at central distances
  $d$ (\ref{eq:2cd}). 
  Top left: $d>r_1+r_2,\,\, r^2<0$. 
  Top center: $r_2<d<r_2+r_1,\,\, r^2>0$.
  Top right: $r_2-r_1<d<r_2,\,\, r^2>0$.
  Bottom left: $d+r_1<r_2,\,\, r^2<0$.
  Bottom center: smaller $d$.
  Bottom right: similar to bottom left, but $C_1$ on other side of $C_2$. }
  \label{fig:2cs}
\end{center}
\end{figure}

We therefore observe (comp. Fig.~\ref{fig:2c}, \ref{fig:2cs} and \ref{fig:r2s}) that (\ref{eq:cchyp})
\begin{itemize}
\item describes for $r^2>0$ all points of real intersection (on the circle $V_1$)
of the two circles.
\item For $r^2<0$ (\ref{eq:cchyp}) becomes the locus equation of the virtual points of 
intersection, i.e. two hyperbola branches that extend symmetrically on both sides of the 
circle $V_1$ (assuming e.g. that we move $V_2$ relative to $V_1$).
\item The sequence of circle meets of Fig.~\ref{fig:2cs} clearly illustrates, that as e.g.
circle $V_2$ moves from the right side closer to circle $V_1$, also the virtual intersection points
approach along the hyperbola branch $(r^2<0)$ on the same side of $V_1$, until the point of 
outer tangence $(d=r_1+r_2,\,\,r=0)$. Then we have real intersection points $(r^2>0)$ until inner tangence
occurs $(d=r_2-r_1, \,\,r=0$). Reducing $d<r_2-r_1$ even further leads to virtual intersection points, 
wandering outwards
on the same side hyperbola branch as before until $d$ becomes infinitely small. Moving $C_2$ over 
to the other side of $C_1$ repeats the phenomenon just described on the other branch of the hyperbola (symmetry
to the vertical symmetry axis of the hyperbola through $C_1$).
\item The transverse symmetry
axis line of the two hyperbola branches is given by $C_1 \wedge C_2\wedge  \bv{n} $, 
i.e. the straight line through the two circle centers. 
\item The assymptotics are at angles $\pm \frac{\pi}{4}$ to the symmetry axis. 
\item The radius $r_1$ is the semitransverse axis segment.
\end{itemize}
%
%


\section{Full Meet of Circle and Straight Line in One Plane}

\begin{figure}
\begin{center}
\includegraphics[width=0.75\textwidth]{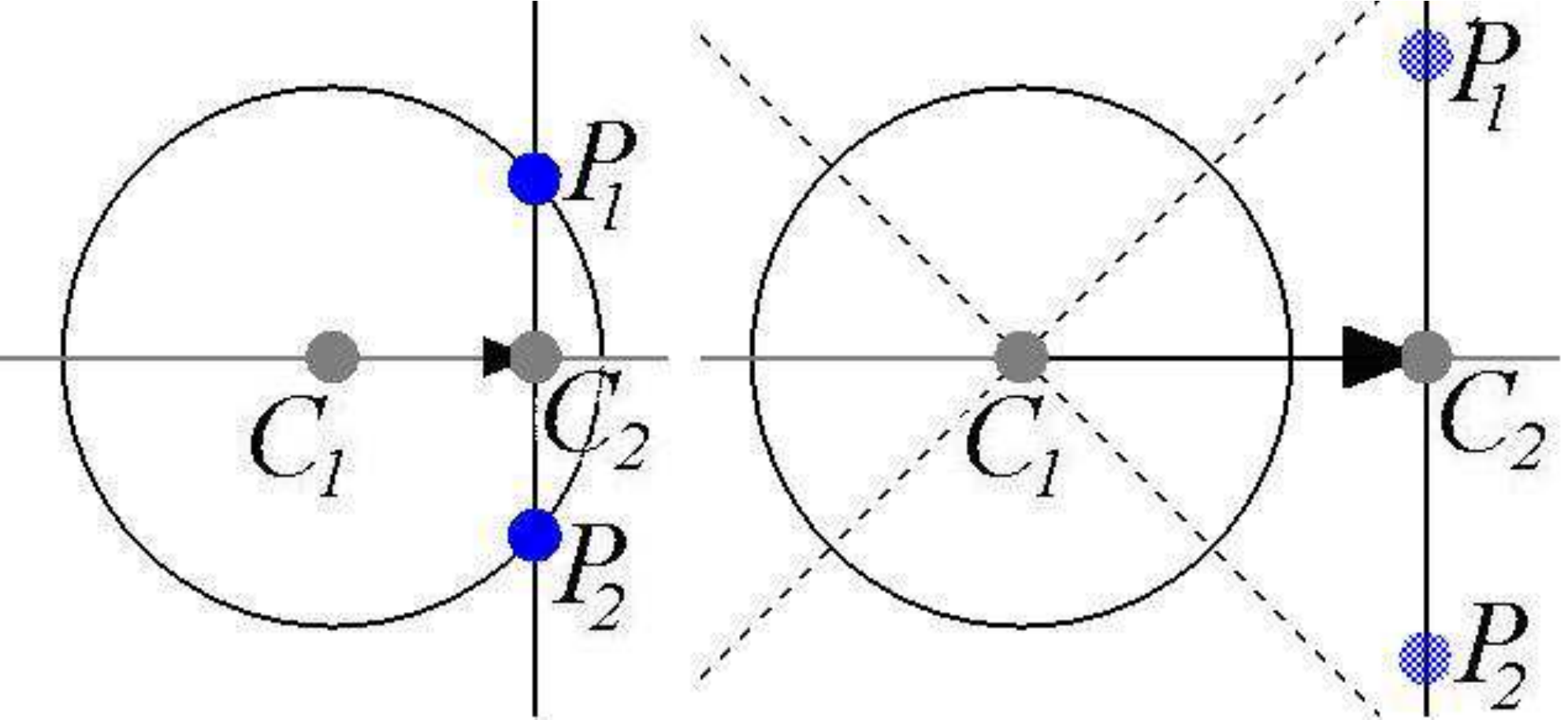}
  \caption[r square sign]{
  Real and virtual intersections of circle and straight line.}
  \label{fig:cln}
\end{center}
\end{figure}

Now we turn our attention to the meet of a
circle with center $C_1$, radius $r_1$ in plane $I_c,\,\, I_c^2=-1$ 
(including the origin $\bv{\bar{n}}$)
\be
  V_1 = (C_1-\frac{1}{2}r_1^2\bv{n})I_c N \enspace,
\ee
and a straight line through $C_2$, with direction $\bv{\hat{p}}$ and in 
the same plane $I_c$,
\be
  V_2 = \bv{\hat{p}}\wedge C_2 \wedge \bv{n} 
      = \bv{\hat{p}}\wedge \bv{c}_2\bv{n}-  \bv{\hat{p}}N \enspace.
\ee
(For convenience $C_2$ be selected
such that $d=\mid\bv{c}_2-\bv{c}_1 \mid$ is the 
distance of the circle center $C_1$ from the line $V_2$. 
See Fig.~\ref{fig:cln}.)

The meet of $V_1$ and $V_2$ is
\be
  M=(V_1\LC J^{-1})\LC V_2 = \ldots = \frac{-1}{2r} P_1\wedge P_2
  \label{eq:clmeet}
\ee
with the join $J=I_c N$,
\be
  r^2 = M^2 =  r_1^2-d^2 \enspace,
  \label{eq:clhyp}
\ee
and
\be
  \bv{p}_1=\bv{c}_2+r\bv{\hat{p}}, \,\,\,
  \bv{p}_2=\bv{c}_2-r\bv{\hat{p}} \enspace.
\ee
$r$ and $\bv{\hat{p}}$ have the same meaning as in (\ref{eq:PPrpc}).
Note that $r^2<0$ for $d>r_1$.

We observe that
\begin{itemize}
\item Equations (\ref{eq:cchyp}) and (\ref{eq:clhyp}) are remarkably similar.
\item The point pair $P_1 \wedge P_2$ is \textbf{now always} on the straight line $V_2$, 
and has center $C_2$!
\item For  $r^2>0$ $(d<r_1)$ 
\be
  r^2+d^2 = r_1^2
\ee
describes the real intersections of circle and straight line.
\item For  $r^2<0$ $(d>r_1)$ 
\be
  r^2+d^2 = r_1^2
\ee
describes the virtual intersections of circle and straight line.
\item The general formula $M\wedge \bv{n} = -V_2$ 
holds for all values of $r$, even if 
$r^2=M^2=0 \,\,(\bv{p}_1=\bv{p}_2=\bv{c}_2)$. 
In this special case $V_2$ is tangent to the circle.\footnote{
For the case of tangency we have now
$\lim_{r \rightarrow 0}M=\,\,-\{\bv{\hat{p}}+C_2\ast \bv{\hat{p}}\,\,\bv{n} \}\wedge C_2 $,
which can be interpreted~\cite{LD:ITMga2003} as a vector in the line $V_2$ attached to $C_2$, 
tangent to the circle.
}
\item In all other repects, the virtual intersection locus hyperbola has the
same properties 
(symmetry, transverse symmetry axis, assymptotics and semitransverse axis segment) 
as that of the meet of two circles in one plane.
\end{itemize}

\section{Full Meet of Two Spheres}

\begin{figure}
\begin{center}
\includegraphics[width=1.0\textwidth]{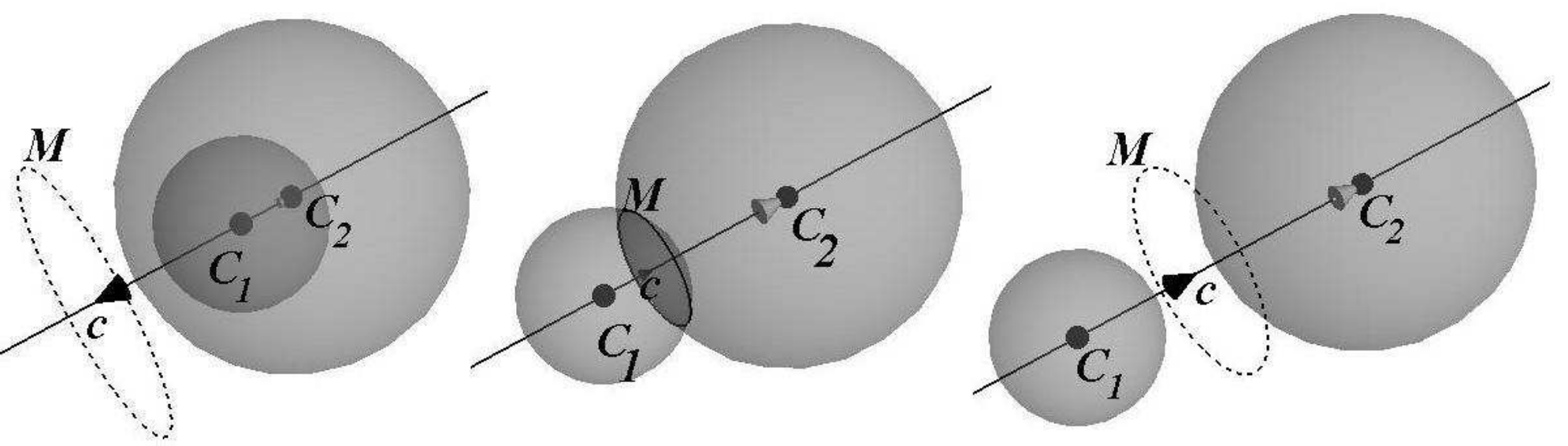}
  \caption[two spheres]{
  Real and virtual intersections of two spheres ($r_1<r_2$).
  Left: $d<r_2-r_1<r_2$, center: $r_2<d<r_1+r_2$, right: $d>r_1+r_2$.}
  \label{fig:2sint}
\end{center}
\end{figure}

Let us assume two spheres (see Fig.~\ref{fig:2sint})
\be
  V_1 = (C_1-\frac{1}{2}r_1^2\bv{n})I N, \,\,\,\,\,\,
  V_2 = (C_2-\frac{1}{2}r_2^2\bv{n})I N \enspace.
  \label{eq:2sV1V2}
\ee
with centers $C_1, C_2$, radii $r_1, r_2$, 
and (3+2)--dimensional pseudoscalar join $J = I N$.
The meet of these two spheres is
$$
  M = (V_1\LC J^{-1})\LC V_2
$$
$$
  =\frac{1}{2}[(\bv{c}_1^2-r_1^2)-(\bv{c}_2^2-r_2^2)]I
   -\frac{1}{2}[(\bv{c}_2^2-r_2^2)\bv{c}_1 - (\bv{c}_1^2-r_1^2)\bv{c}_2 ]I \bv{n}
   +(\bv{c}_1-\bv{c}_2)I\bv{\bar{n}}
   +\bv{c}_1\wedge \bv{c}_2 I N 
$$
\be
  = \ldots = 
  \,\,d (C+\frac{1}{2}r^2\bv{n})\wedge \{I_c+\bv{n}(C\LC I_c)  \} 
  \enspace,
  \label{eq:2smeet}
\ee
where $r$ and $d$ are defined as for the case of intersecting two circles
[$r^2 = -M^2/(M\wedge \bv{n})^2$, note the sign!]
We further introduced in (\ref{eq:2smeet}) the plane bivector 
\be
  I_c = \frac{(\bv{c}_1-\bv{c}_2)}{d}I
  \label{eq:2sIc}
\ee
and the vector (see Fig.~\ref{fig:2sint})
\be
  \bv{c}= \bv{c}_1+\frac{1}{2}(1+\frac{r_1^2-r_2^2}{d^2})(\bv{c}_2-\bv{c}_1)
  \enspace.
\ee
Comparing (\ref{eq:cfcirc}) and (\ref{eq:2smeet}) we see that 
$M$ is a conformal circle multivector with 
radius $r$,
oriented parallel to $I_c$ in the plane 
$M\wedge \bv{n}=d C \wedge I_c\wedge \bv{n}$,
and with center $C$.

\begin{figure}
\begin{center}
\includegraphics[width=0.4\textwidth]{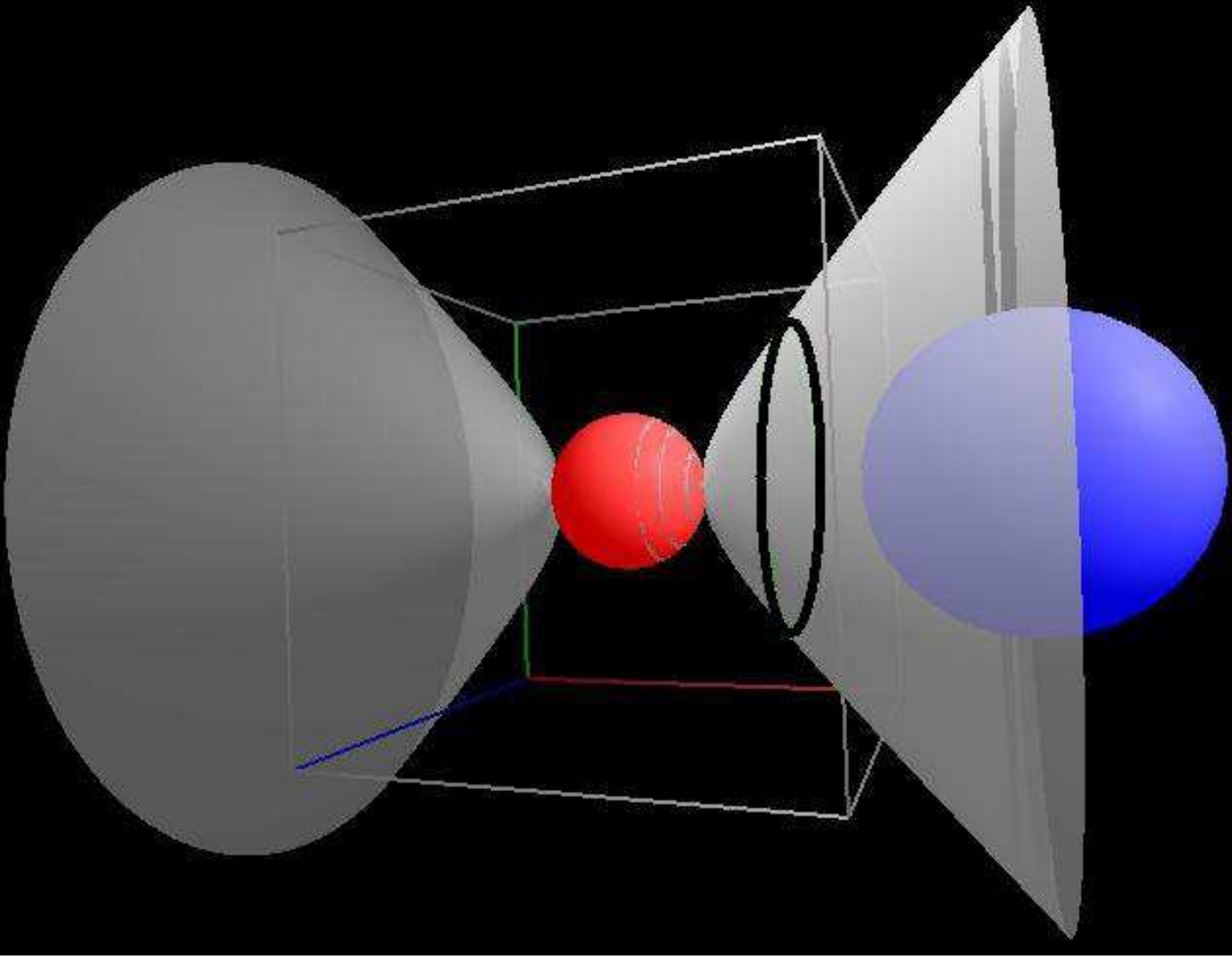}
  \caption[hyperboloid]{
  Two sheet hyperboloid of virtual intersections of two spheres.}
  \label{fig:hypv}
\end{center}
\end{figure}

Regarding the formula
\be
  r^2+d_1^2=r_1^2
  \label{eq:2sr}
\ee
with $d_1=|\bv{c}-\bv{c}_1 | $ it remains to observe that
\begin{itemize}
\item equation (\ref{eq:2sr}) describes for $r^2>0$ $(d_1<r_1)$ the real radius $r$ 
circles of intersection of two spheres.
\item These intersection circles are centered at 
$
  C=\bv{c}+\frac{1}{2}c^2\bv{n} + \bv{\bar{n}}
$
in the plane $M\wedge \bv{n}$ perpendicular to the center connecting straight line 
$C_1\wedge C_2 \wedge \bv{n}$, i.e. parallel to the bivector of (\ref{eq:2sIc}).
\item For $r^2=0$, $M\wedge \bv{n}$ gives still the (conformal) tangent plane trivector 
of the two spheres.\footnote{
For the case of tangency $(r^2=0)$ we have now
$M = d \,\,C\wedge \{I_c+\bv{n}(C \LC I_c)  \}$, which can be interpreted~\cite{LD:ITMga2003} 
as tangent direction bivector $I_c$ of the two
tangent spheres, located at the point $C$ of tangency.
}
\item We have for $r^2<0$ $(d_1>r_1)$ virtual circles of intersection 
forming a hyperboloid with \textbf{two sheets}, as shown in Fig.~\ref{fig:hypv}.
\item The transverse symmetry axis (straight) line of the two sheet hyperboloid 
is $C_1\wedge C_2 \wedge \bv{n}$.
\item The discussion of the meet of two circles in section \ref{sc:2C} 
related to Fig.~\ref{fig:2cs} and Fig.~\ref{fig:r2s}
applies also to the case of the meet of two spheres. $r^2$ is now the squared 
radius of real and virtual meets (circles instead of point pairs). 
\item The asymptotic double cone has angle $\pi/4$ relative to the transverse
symmetry axis.
\item The sphere radius (e.g. $r_1$) is again the semitransverse axis segment
 of the two sheet hyperboloid
(assuming e.g. that we move $V_2$ relative to $V_1$). 
\end{itemize}

\section{Full Meet of Sphere and Plane}

\begin{figure}
\begin{center}
\includegraphics[width=0.75\textwidth]{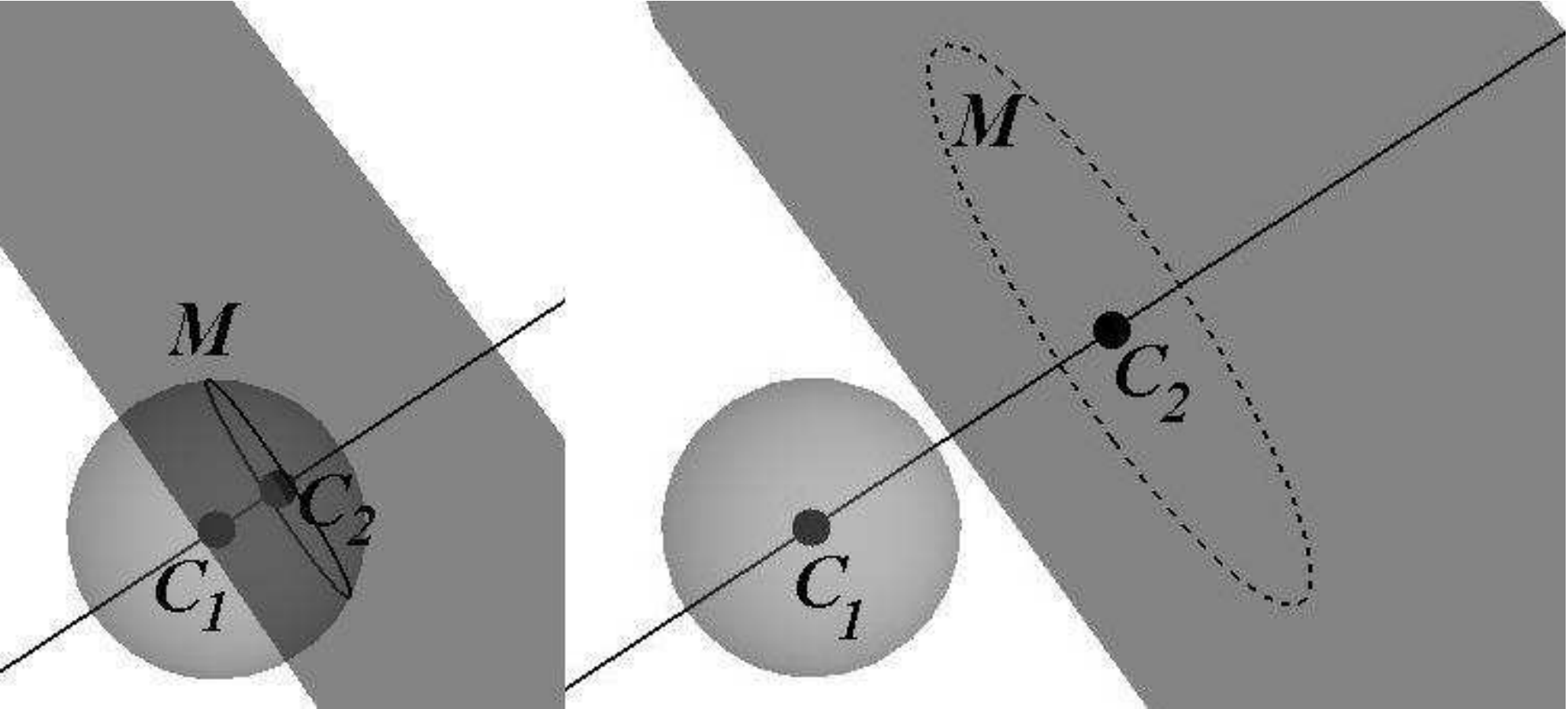}
  \caption[two spheres]{
  Real and virtual intersections of sphere and plane.}
  \label{fig:2sp}
\end{center}
\end{figure}

Let $V_1$ be a conformal sphere four-vector [as in (\ref{eq:2sV1V2})] 
and $V_2$ a conformal plane four-vector 
[(\ref{eq:cplane}) with normalization $\alpha=1$]. 
The analogy to the case of circle and straight line is now obvious. 
For the virtual $(r^2<0)$ intersections we get the same 
two sheet hyperboloid as for the case of sphere and sphere,
but now both real and virtual intersection circles are always on the plane $V_2$,
as shown in Fig.~\ref{fig:2sp}.
The general formula $M\wedge\bv{n}= V_{2}$ 
holds for all values of $r$, even if $r^2=M^2=0$. 
In this special case $V_2$ is tangent to the sphere.

\section*{Acknowledgement}

Soli Deo Gloria. I thank my wife and my children. I further
thank especially Prof. Hongbo Li and his colleagues for organizing the GIAE workshop. 
The GAViewer~\cite{LD:GAV} was used to create 
Figs.~\ref{fig:cpp},\ref{fig:2c},\ref{fig:2cs},\ref{fig:cln},\ref{fig:2sint} and \ref{fig:2sp} 
and to probe many of the formulas. 
In Fig.~\ref{fig:r2s} the interactive geometry software Cinderella~\cite{RGK:cindy} was used. 
Fig.~\ref{fig:hypv} was created by C. Perwass with
CLUCalc~\cite{CP:CLU}. I thank the Signal Processing Group of the CUED for its hospitality
in the period of finishing this paper and L. Dorst for a number of important comments.

\end{document}